\documentclass[12pt]{amsart}
\usepackage{amssymb}

\newtheorem{theorem}{Theorem}[section]

\numberwithin{equation}{section}

\begin{document}

\baselineskip=15pt
 
\title[Vector bundles on Symmetric product of curves.]{Vector bundles on
Symmetric product of curves.}

\author[D. S. Nagaraj]{D. S. Nagaraj}

\address{The Institute of Mathematical Sciences, CIT
Campus, Taramani, Chennai 600113, India}

\email{dsn@imsc.res.in}



\date{}

\maketitle

\section{Introduction}

 The aim of the article  is to give a survey of some recent work on the
properties some special vector bundles on symmetric product of curves.
The title and content of the article is almost same as that of the my talk
at the "Manipal Workshop on Algebraic Geometry", held at Manipal during
5-9th Jan 2015. The references are 
added so that interested reader can find the details about the 
results mentioned here. 
\section{symmetric product of curves}
Let $C$ be a smooth irreducible curve over the field $\mathbb{C}$ of complex
numbers. Let $n$ be an integer greater or equal to $2.$ Then the symmetric group
$S_n$ acts naturally on $C^n:= C\times \cdots \times C$ - the n-fold product of
$C.$ Then  
the quotient space is a smooth projective variety of dimension $n$ and is
denoted by $C_n.$ A point of $C_n$ can be identified with a effective divisor
on $C$ of degree $C$ and if it is the image of $(P_1,\ldots, P_n)$ and is 
denoted by $P_1+\cdots + P_n.$ Thus 
$C_n$ can be identified with set of all
effective divisors of degree $n$ on $C.$ The subset $\{P_1,\ldots, P_n\}$
of $C$ is called the support of the effective divisor $P_1+\cdots + P_n$
of degree $n.$

For example when $C = \mathbb{P}^1$ the projective line the projective
manifold $C_n$ can be seen to be isomorphic to $\mathbb{P}^n$ the projective 
$n$ space over the field of complex numbers.
If $C$ is an elliptic curve over $\mathbb{C},$ i.e., a curve of genus one over
$\mathbb{C}$ then using the group law we get a morphism
$C_n \to C$ sending the effective divisor $D= P_1+\cdots P_n$ to the
corresponding
sum under the group law of the curve $C.$ Then one get by Riemann-Roch theorem
$C_n$ as a $\mathbb{P}^{n-1}$ bundle over $C.$

Let $C_n$ be the $n$th symmetric product of smooth irreducible curve $C.$
Consider the subset $\Delta_n \subset C_n\times C$ consisting of all
$(D, P)\in C_n\times C$ where $D$ thought of an effective divisor of degree 
$n$ on $C$ contains the point $P$ in its support. It is easy to see that
$\Delta_n$ has natural structure of a reduced divisor on $\Delta_n\times
C.$ It is easy to see that $\Delta_n$ is a smooth divisor. In fact one 
can show that $\Delta_n$ is isomorphic to $C_{n-1}\times C.$  
This divisor is called universal effective divisor of degree $n.$ The projection 
$p_1: C_n\times C \to C_n$ restricted to $\Delta_n$ gives a finite 
morphism $q: \Delta_n \to C_n$ of degree $n.$
For each $x\in C$ there is a natural divisor 
$$x+C_{n-1}=\{x+x_1+\cdots+x_{n-1}| x_1+\cdots+x_{n-1}\in C_{n-1} \}.$$
The inverse image of $x+C_{n-1}$ under the natural quotient map
$C^n \to C_n$ is equal to 
$$\sum_{i=1}^n p_i^*(x), $$
where $p_i:C^n \to C $ is the $i$-th projection ($1\leq i \leq n$).
Since $\sum_{i=1}^n p_i^*(x) $ is an ample divisor on $C^n$ we see that
$x+C_{n-1}$ is an ample divisor on $C_n.$

\section{Special bundles on $C_n$ and their properties}

Let $E$ be any vector bundle of rank $r$ on $C.$ Then set
$$\mathcal{F}_n(E)= q_*(p_2^*(E)|_{\Delta_n}),$$
where $p_2: C_n\times C \to C$ is the second projection and 
$q: \Delta_n \to C_n$ is the restriction of the  first projection 
$p_1: C_n\times C \to C_n.$ Since $q: \Delta_n \to C_n$ is a
finite morphism of degree $n,$ the sheaf $\mathcal{F}(E)$ is a
vector bundle of rank $nr.$ 

A study of the secant vector bundles on $C_2$ was initiated
by R. L. E. Schwarzenberger in \cite{Sc}.
In \cite{BL} the bundle $\mathcal{F}_n(E)$ is considered and 
studied its properties as a natural parabolic bundle.
If  $L$ is a line bundle on $C$ the bundle bundle $\mathcal{F}_n(L)$ 
on $C_n$ defined above  is known as the $n$-th secant bundle corresponding
to $L.$ For more details about these bundles we refer to \cite[chapter VIII]{ACGH}. 

It is natural to ask what properties of
$E$ will be inherited  by $\mathcal{F}_n(E).$ In particular 
we can ask the following questions:
\begin{itemize}
\item[1)] If $E$ is stable (respectively, semi-stable) on $C,$ then does 
it imply $\mathcal{F}_n(E)$
is stable (respectively, semi-stable) on $C_n?$ 
\item[2)] If $E$ and $F$ are two bundles on $C$ such that 
$\mathcal{F}_n(E)\simeq\mathcal{F}_n(F)$ on $C_n,$ then does it 
imply $E \simeq F$?
\end{itemize}

Recall that for a vector bundle $V$ on a curve $X$ is stable
(respectively, semi-stable) if for any sub bundle  $W$ of $V$
$$\mu(W)(:=\frac{{\rm degree}(W)}{{\rm rank}(W)}) < \mu(V)$$
(respectively, $\mu(W)\leq \mu(V)$).
Let $V$ be  torsion free sheaf  on smooth projective variety $X$ of
dimension $n.$ Let $[H]$ be the cohomology class
of an ample line bundle $H$ on $X.$ Recall that the degree of $V$ 
with respect to $H,$ denoted by ${\rm degree}_H(V),$ 
is defined to be the number $C_1(V).[H]^{n-1},$ where $C_1(V)$ is the first
Chern class of $V$ and 'dot' denote the cup product in the cohomology
ring. The slope $\mu_H(V)$ of $V$ with respect to $H$ is defined to be 
${\rm degree}_H(V)/{{\rm rank}(V)}.$
For a vector bundle (more generally a torsion free sheaf) on $X$
is said to be stable (respectively, semi-stable) with respect to $H,$ if for any  
$\mathcal{O}_X$ sub sheaf  $W$ of $V$
$$\mu_H(W) < \mu_H(V)$$
(respectively, $\mu_H(W)\leq \mu_H(V)$).

In question 1) above  we take the ample line bundle $H$ to be the line bundle 
associated to the ample divisor $x+C_{n-1}.$

In general both the questions have negative answer.
 
 We give an example to show that in general 2) does not always hold. 
 
 We show that there vector bundles $E, F$ on $\mathbb{P}^1$ such that $\mathcal{F}_2(E)\,=\,
 \mathcal{F}_2(F)$ on $(\mathbb{P}^1)_2\simeq \mathbb{P}^2$ but
imply that $E\,=\, F$.

Note that $(\mathbb{P}^1)_2\,\simeq \,\mathbb{P}^2.$
If we identify $(\mathbb{P}^1)_2$ with $\mathbb{P}^2,$ then the universal degree two divisor
$$
{\Delta}_2\, \subset\, \mathbb{P}^1\times (\mathbb{P}^1)_2
\,\simeq\, \mathbb{P}^1\times \mathbb{P}^2
$$
is the zero locus of a section of the line bundle $p^*(\mathcal{O}_{\mathbb{P}^1}(2))\otimes
q^*(\mathcal{O}_{\mathbb{P}^2}(1)),$ where
\begin{equation}\label{ee1}
p\, :\, \mathbb{P}^1\times \mathbb{P}^2
\, \longrightarrow\, \mathbb{P}^1\ ~ \text{ and } ~ \ q \, :\,
\mathbb{P}^1\times \mathbb{P}^2
\,\longrightarrow\, \mathbb{P}^2
\end{equation}
are the natural projections. From  this we see that
\begin{itemize}
\item $\mathcal{F}_2(\mathcal{O}_{\mathbb{P}^1}(1))
\,=\, \mathcal{O}_{\mathbb{P}^2}\oplus \mathcal{O}_{\mathbb{P}^2}$

\item $\mathcal{F}_2(\mathcal{O}_{\mathbb{P}^1}(-1))
\,=\, \mathcal{O}_{\mathbb{P}^2}(-1)\oplus \mathcal{O}_{\mathbb{P}^2}(-1)$

\item $\mathcal{F}_2(\mathcal{O}_{\mathbb{P}^1})
\,=\, \mathcal{O}_{\mathbb{P}^2}\oplus \mathcal{O}_{\mathbb{P}^2}(-1)$.
\end{itemize}
For any two vector bundles $E$ and $F$ on
$\mathbb{P}^1$ we have 
$\mathcal{F}_2(E\oplus F) \,=\, \mathcal{F}_2(E) \oplus \mathcal{F}_2(F)$ on $\mathbb{P}^2.$ 
From these observations it follows that
$$\mathcal{F}_2(\mathcal{O}^{\oplus 2m}_{{\mathbb P}^1})
\,=\, \mathcal{O}_{\mathbb{P}^2}^{\oplus 2m}\oplus 
\mathcal{O}_{\mathbb{P}^2}(-1)^{\oplus 2m}
\,=\,\mathcal{F}_2(\mathcal{O}_{\mathbb{P}^1}(1)^m\oplus \mathcal{O}_{\mathbb{P}^1}(-1)^m)\, .$$

To give negative answer to the question 1) we note that any line bundle
on a curve is stable and direct sum of copies of a same line bundle
is semistable. With above notation we see that 
$$\mathcal{F}_2(\mathcal{O}_{\mathbb{P}^1}) = \mathcal{O}_{\mathbb{P}^2}
\oplus \mathcal{O}_{\mathbb{P}^2}(-1)$$
and 
$$\mathcal{F}_2(\mathcal{O}^m_{\mathbb{P}^1}) = \mathcal{O}^m_{\mathbb{P}^2}
\oplus \mathcal{O}_{\mathbb{P}^2}(-1)^m$$
and the bundle $\mathcal{O}^m_{\mathbb{P}^2}
\oplus \mathcal{O}_{\mathbb{P}^2}(-1)^m$ is not semistable for all $m\geq 1.$

The question 1) is treated in several papers, for
example see \cite{KS} \cite{BN1},\cite{ELN}, \cite{GuHe}. 

One of best known result 
in this direction is the following:

\begin{theorem}\label{thm1}
Let $L$ be any nontrivial line bundle on a smooth connected complex
projective curve $C.$
\begin{enumerate}
\item The vector bundle $\mathcal{F}_2(L)$ is semistable on $C_2.$ 

\item The vector bundle $\mathcal{F}_2(L)$ is stable 
 $C_2,$ unless $L$ is of the form ${\mathcal O}_C(x)$
or ${\mathcal O}_C(-x)$ for some point $x\,\in\, C$.
\end{enumerate}
\end{theorem}

{\bf Proof} See \cite[Theorem 3.1.]{BN1}.

The stability of the bundles $\mathcal{F}_2(\mathcal{O}_{\mathbb{P}^1}(m))$ 
for $m\geq 2$ was proved in \cite{GuHe}. The above theorem with
more assumptions was proved in \cite{ELN}.

The question 2) is treated in several papers for example \cite{BN1}
\cite{BN2}, and \cite{BP}.

The best known result about question 2) is the following:
\begin{theorem}\label{thm2}
Let $E_1$ and $E_2$ be stable vector bundles on $C.$ 
Assume $g(C)\geq 2.$ If
$\mathcal{F}_n(E_1)\simeq \mathcal{F}_n(E_2)$
on $C_n$ then $E_1 \simeq E_2$ on $C.$
\end{theorem}

{\bf Proof} See \cite[Proposition 3.1.]{BN2}.

In \cite{BN1}, it is proved  that if $L, L'$ be two line bundles
on a smooth curve $C$ of genus $g\geq 0$ such that
$\mathcal{F}_2(L)\simeq \mathcal{F}_2(L')$
on $C_2$ then $L \simeq L'$ on $C.$  Theorem \ref{thm2} was proved
under more restrictive assumptions in \cite{BP}.

In a recent preprint \cite{BN3}  Theorem \ref{thm2} was improved
in the sense that Theorem holds for bundles with bounded "Harder-Narasimhan"
filtration depending on the genus of the curve. This result is used to prove the following:
  
\begin{theorem}\label{thm3}
Let $E_1$ and $E_2$ be  vector bundles on a smooth projective surface $S.$ 
 If
$\mathcal{H}_n(E_1)\simeq \mathcal{H}_n(E_2)$
on ${\rm Hilb}_n(S)$ then $E_1 \simeq E_2$ on $S.$
\end{theorem}

Here ${\rm Hilb}_n(S)$ denote the Hilbert Scheme that parameterizes the zero dimensional  
closed sub schemes of length $n$ on $S$ and $\mathcal{H}_n(E)$ is the vector
bundle on ${\rm Hilb}_n(S)$ associated to a vector bundle $E$ on $S.$
The construction of the bundle $\mathcal{H}_n(E)$ is similar to the
above construction of bundle  $\mathcal{F}_n(V).$

\end{document}